\title{Eigenvalues and Homology of Flag Complexes \\
and Vector Representations of Graphs}
\documentstyle[amsmath,amsfonts,amssymb,12pt]{article}

\font\sixbb=msbm6
\font\eightbb=msbm8
\font\twelvebb=msbm10 scaled 1095
\newfam\bbfam
\textfont\bbfam=\twelvebb \scriptfont\bbfam=\eightbb
                           \scriptscriptfont\bbfam=\sixbb
\def\bb{\fam\bbfam\twelvebb}
\newcommand{\Rea}{{\bb R}}

\newcommand{\Rat}{{\bb Q}}
\newcommand{\Int}{{\bb Z}}

\newtheorem{theorem}{\bf Theorem}[section]
\newtheorem{claim}[theorem]{\bf Claim}

\newtheorem{proposition}[theorem]{\bf Proposition}

\newcommand{\enp}{\begin{flushright} $\Box$ \end{flushright}}
\newcommand{\beq}[0]{\begin{equation}}
\newcommand{\enq}[0]{\end{equation}}
\newcommand{\phu}{\phi_u}

\newcommand{\lk}{{\rm lk}}
\newcommand{\th}{{\tilde H}}

\newcommand{\cf}{{\cal F}}

\newcommand{\dkos}{d_{k-1}^*}
\newcommand{\comg}{\overline{G}}

\newcommand{\vone}{{\bf 1}}

\newcommand{\valpha}{{\bf \alpha}}

\newcommand{\cfi}{\{\cf_i\}_{i=1}^m}

\newcommand{\ci}{{\rm I}}
\newcommand{\cig}{\ci(G)}

\newcommand{\cm}{\mathcal M}

\newcommand{\tg}{\tilde{\gamma}}

\newcommand{\va}{{\bf a}}

\newcommand{\gva}{G_{\va}}

\newcommand{\ig}{i\gamma}

\begin{document}
\author{{\sc R. Aharoni}$^*$ \and {\sc E. Berger}$^*$ \and
{\sc R. Meshulam}\thanks{Department of Mathematics, Technion,
Haifa 32000, Israel.~~  e-mails:
\newline  ra@tx.technion.ac.il,
eberger@princeton.edu, meshulam@math.technion.ac.il} }
\date{}
\insert\footins{\footnotesize\rule{0pt}{\footnotesep}
\\ {\it Math Subject Classification:}  13F55, 05C69.
}
\maketitle
\pagestyle{plain}
\begin{abstract}
The flag complex of a graph $G=(V,E)$ is the simplicial complex
$X(G)$ on the vertex set $V$ whose simplices are subsets of $V$
which span complete subgraphs of $G$. We study relations between
the first eigenvalues of successive higher Laplacians of $X(G)$.
One consequence is the following \ \\ \\{\bf Theorem:} Let
$\lambda_2(G)$ denote the second smallest eigenvalue of the
Laplacian of $G$. If $\lambda_2(G)
> \frac{k}{k+1} |V|$ then $\th^{k}(X(G);\Rea)=0$. \ \\ \\
Applications include a lower bound on the homological connectivity
of the independent sets complex $\ci(G)$, in terms of a new graph
domination parameter $\Gamma(G)$ defined via certain vector
representations of $G$. This in turns implies a Hall type theorem
for systems of disjoint representatives in hypergraphs.

\end{abstract}

\section{Introduction}

Let $G=(V,E)$ be a graph with $|V|=n$ vertices. The {\it
Laplacian} of $G$ is the $V \times V$ positive semidefinite matrix
$L_G$ given by $$L_G(u,v) = \left\{
\begin{array}{ll}
        \deg(u)  & u=v \\
        -1       & uv \in E \\
         0   & {\rm otherwise}
\end{array}
\right.~~ $$ Let $0=\lambda_1(G) \leq \lambda_2(G) \leq \cdots
\leq \lambda_n(G)$ denote the eigenvalues of $L_G$. The second
smallest eigenvalue $\lambda_2(G)$, called the {\it spectral gap},
is a parameter of central importance in a variety of problems. In
particular it controls the expansion properties of $G$ and the
convergence rate of a random walk on $G$ (see e.g. \cite{Boll98}).
The {\it Flag Complex} of $G$ is the simplicial complex $X(G)$ on
the vertex set $V$ whose simplices are all subsets $\sigma \subset
V$ which form a complete subgraph of $G$. Topological properties
of $X(G)$ play key roles in recent results in matching theory (see
below).
\\
In this paper we study relations between $\lambda_2(G)$, the
cohomology of $X(G)$, and a new graph domination parameter
$\Gamma(G)$ which is defined via certain vector representations of
$G$. As an application we obtain a Hall type theorem for systems
of disjoint representatives in  families of hypergraphs. \ \\ \\
For $k \geq -1$ let $C^k(X(G))$ denote the space of real valued
simplicial $k$-cochains of $X(G)$ and let $d_k:C^k(X(G))
\rightarrow C^{k+1}(X(G))$ denote the coboundary operator. For $k
\geq 0$ define the reduced $k$-dimensional Laplacian of $X(G)$ by
$\Delta_k=d_{k-1}d_{k-1}^*+d_k^*d_k$ (see section \ref{hodge} for
details). Let $\mu_k(G)$ denote the minimal eigenvalue of
$\Delta_k$. Note that $\mu_0(G)=\lambda_2(G)$.
Our main result is
the following
\begin{theorem}
\label{eigenv} For $k \geq 1$
\begin{equation}
\label{maineq} k\mu_k(G) \geq (k+1)\mu_{k-1}(G) -n~~.
\end{equation}
\end{theorem}
As a direct consequence of Theorem \ref{eigenv} we
obtain
\begin{theorem} \label{ei} If $\lambda_2(G) >
\frac{kn}{k+1}$ then $\th^{k}(X(G),\Rea)=0$.
\end{theorem}
{\bf Remarks:} 
\vspace{5pt}
\\
1. Theorem \ref{ei} is related to a well-known result
 of Garland 
(Theorem 5.9 in \cite{G73}) and its extended version
by Ballmann and \'{S}wi\c{a}tkowski
(Theorem 2.5 in \cite{BS97}). 
Roughly speaking, these results (in their simplest untwisted
form)  guarantee the vanishing of $\th^k(X;\Rea)$ 
provided that for {\it each} $(k-1)$-simplex $\tau$ in $X$, 
the spectral gap of the $1$-skeleton of the link of $\tau$
is sufficiently large. Theorem \ref{ei} is, in a sense, a global
counterpart of this statement for flag complexes.
\vspace{5pt}
\\
2. Let $n=r \ell$ where $r \geq 1,\ell \geq 2$, and let $G$
be the Tur\'an graph $T_r(n)$, i.e. the complete $r$-partite graph
on $n$ vertices with all sides equal to $\ell$. The flag complex
$X(T_r(n))$ is homotopy equivalent to the wedge of $(\ell-1)^r$
$(r-1)$-dimensional spheres. It can be checked that
$\mu_k(T_r(n))=\ell(r-k-1)$ for all $0 \leq k \leq r-1$, hence
 (\ref{maineq}) is satisfied with equality.
Furthermore, $\lambda_2(G)=\ell(r-1)=\frac{r-1}{r}n$ while
$\th^{r-1}(X(G)) \neq 0$. Therefore the assumption in Theorem
\ref{ei} cannot be replaced by $\lambda_2(G) \geq \frac{kn}{k+1}$.
\ \\ \\ We next study some graph theoretical consequences of
Theorem \ref{ei}. The {\it Independence Complex} $\ci(G)$ of $G$
is the simplicial complex on the vertex set $V$ whose simplices
are all independent sets $\sigma \subset V$. Thus
$\ci(G)=X(\comg)$ where $\comg$ denotes the complement of $G$.
Recent work on hypergraph matching, starting in \cite{AH00} with
later developments in \cite{Aharoni01,M01,ABZ02,ACK02,M03}, has
utilized topological properties of $\ci(G)$ to derive new Hall
type theorems for hypergraphs. The main ingredient in these
developments are lower bounds on the homological connectivity of
$\ci(G)$. For a simplicial complex $Z$ let
$\eta(Z)=\min\{i:\th^i(Z,\Rea) \neq 0\}+1~.$ It turns out that
various domination parameters of $G$ may be used to provide lower
bounds on $\eta(\ci(G))$. For a subset of vertices $S \subset V$
let $N(S)$ denote all vertices that are adjacent to at least one
vertex of $S$ and let $N'(S)=S \cup N(S)$. $S$ is a {\it
dominating set} if $N'(S) = V$. $S$ is a {\it totally dominating
set} if $N(S) = V$. Here are a few domination parameters:
\begin{itemize}
\item
The {\it domination number} $\gamma(G)$ is the minimal size of a
dominating set.
\item
The {\it total domination number} $\tg(G)$ is the minimal size of
a totally dominating set.
\item
The {\it independent domination number} $\ig(G)$ is the maximum,
over all independent sets $I$ in $G$, of the minimal size of a set
$S$ such that $N(S) \supset I$.
\item The {\it strong fractional domination number},
$\gamma^*_s(G)$ is the minimum of $\sum_{v\in V}f(v)$, over all
nonnegative functions $f:V \to \mathbb{R}$ such that \\ $\sum_{uv \in
E}f(u)+\deg(v)f(v) \ge 1$ for every vertex $v$.
\end{itemize}
Some known lower bounds on $\eta$ are: $\eta(\cig) \ge \tg(G)/2$
\cite{M01}, $\eta(\cig) \ge \ig(G)$ \cite{AH00}, $\eta(\cig) \ge
\gamma^*_s(G)$ \cite{M03}. \ \\ \\ Here we introduce a new
domination parameter, defined by vector representations. It is
similar in spirit to the $\Theta$ function defined by Lov\'asz
\cite{Lovasz79}. It uses vectors to mimick domination, in a
 way similar to that in which the $\Theta$ function mimicks independence of sets
of vertices. It is defined as follows. A {\em vector
representation} of a graph $G=(V,E)$ is an assignment $P$ of a
vector $P(v) \in \Rea^\ell$ for some fixed $\ell$ to every vertex
$v$ of the graph, such that the inner product $P(u) \cdot P(v) \ge
1$ whenever $u,v$ are adjacent in $G$ and $P(u) \cdot P(v) \geq 0$
if they are not adjacent. We shall identify the representation
with the matrix $P$ whose $v$-th row is the vector $P(v)$.
\\
Let $\vone$ denote the all $1$ vector in $\Rea^V$. A non-negative
vector $\valpha$ on $V$ is said to be {\em dominating for $P$} if
$\sum_{v \in V} \alpha(v) P(v) \cdot P(u) \ge 1$ for every vertex
$u$, namely $\valpha PP^T \ge \vone$. (Note that taking $\valpha$
to be the characteristic function of some totally dominating set
satisfies this condition regardless of the representation.) The
{\it value} of $P$ is $$|P|=\min\{ \valpha \cdot \vone~:~\valpha
\geq 0~,~\valpha PP^T \ge \vone~\}~~.$$ The supremum of $|P|$ over
all vector representations $P$ of $G$ is denoted by $\Gamma(G)$.
Our main application of Theorem \ref{ei} is the following
\begin{theorem}\label{etagegamma}
$ \eta(\cig) \ge \Gamma(G)~.$
\end{theorem}
{\bf Remark:} One natural vector representation of $G$ is obtained
by taking $P(v) \in \Rea^E$ to be the edge incidence vector of the
vertex $v$. For this representation $|P|=\gamma^*_s(G)$ hence
$\Gamma(G) \geq \gamma^*_s(G)$. The bound $\eta(\cig) \geq
\gamma^*_s(G)$ was previously obtained in \cite{M03}.  Theorem
\ref{etagegamma} is however stronger and often gives much sharper
estimates for $\eta(\cig)$, see e.g. the case of cycles described
in Section \ref{s:etag}. \ \\ \\ We next use Theorem
\ref{etagegamma} to derive a new Hall type result for hypergraphs.
Let $\cf \subset 2^V$ be a hypergraph on a finite ground set $V$.
The {\it width} $w(\cf)$ of $\cf$ is the minimal $t$ for which
there exist $F_1,\ldots,F_t \in \cf$ such that for any $F \in
\cf$, $F_i \cap F \neq \emptyset$ for some $1 \leq i \leq t$.
\\
The {\it fractional width} $w^*(\cf)$ of $\cf$ is the minimum of
$\sum_{E \in \cf} f(E)$ over all non-negative functions $f : \cf
\rightarrow \mathbb{R}$ with the property that for every edge $E
\in \cf$ the sum $\sum_{F \in \cf} f(F) |E \cap F|$ is at least 1.
A {\it matching} in $\cf$ is a subhypergraph $\cm \subset \cf$
such that $F \cap F' = \emptyset$ for all $F \neq F' \in \cm$. Let
$\cfi$ be a family of hypergraphs. A {\it system of disjoint
representatives (SDR)} of $\cfi$ is a matching $F_1,\ldots,F_m$
such that $F_i \in \cf_i$ for $1 \leq i \leq m$. Haxell
\cite{Haxell95} proved the following
\begin{theorem}{\cite{Haxell95}}
\label{hax95} If $\cfi$ satisfies $w(\cup_{i\in I} \cf_i) \geq
2|I|-1$ for all $\emptyset \neq I \subset [m]$, then $\cfi$ has an
SDR.
\end{theorem}
Here we use Theorem \ref{etagegamma} to show
\begin{theorem}
\label{wstar} If $\cfi$ satisfies $w^*(\cup_{i\in I} \cf_i) >
|I|-1$ for all $\emptyset \neq I \subset [m]$, then $\cfi$ has an
SDR.
\end{theorem}
The paper is organized as follows. In section \ref{hodge} we
recall some topological terminology and the simplicial Hodge
theorem. Theorems \ref{eigenv} and \ref{ei} are proved in section
\ref{s:eandc}. The proofs utilize the approach of Garland
 \cite{G73} and its exposition by
Ballmann and \'{S}wi\c{a}tkowski
\cite{BS97}. In section \ref{s:etag} we relate the $\Gamma$
parameter to homological connectivity and prove Theorem
\ref{etagegamma}. In section \ref{s:hall} we recall a homological
Hall type condition (Proposition \ref{hom}) for the existence of
colorful simplices in a colored complex.  Combining this condition
with Theorem \ref{etagegamma} then completes the proof of Theorem
\ref{wstar}.

\section{Topological Preliminaries}
\label{hodge} Let $X$ be a finite simplicial complex on the vertex
set $V$. Let $X(k)$ denote the set of $k$-dimensional simplices in
$X$, each taken with an arbitrary but fixed orientation. A
simplicial $k$-cochain is a real valued skew-symmetric function on
all ordered $k$-simplices of $X$. For $k \geq 0$ let $C^k(X)$
denote the space of $k$-cochains on $X$. The $i$-face of an
ordered $(k+1)$-simplex $\sigma=[v_0,\ldots,v_{k+1}]$ is the
ordered $k$-simplex
$\sigma_i=[v_0,\ldots,\widehat{v_i},\ldots,v_{k+1}]$. The
coboundary operator $d_k:C^k(X) \rightarrow C^{k+1}(X)$ is given
by $$d_k \phi (\sigma)=\sum_{i=0}^{k+1} (-1)^i \phi
(\sigma_i)~~.$$ It will be convenient to augment the cochain
complex $\{C^i(X)\}_{i=0}^{\infty}$ with the $(-1)$-degree term
$C^{-1}(X)=\Rea$ with the coboundary map $d_{-1}:C^{-1}(X)
\rightarrow C^0(X)$ given by $d_{-1}(a)(v)=a$ for $a \in \Rea~,~v
\in V$. Let $Z^k(X)= \ker (d_k)$ denote the space of $k$-cocycles
and let $B^k(X)={\rm Im}(d_{k-1})$ denote the space of
$k$-coboundaries. For $k \geq 0$ let $\th^k(X)=Z^k(X)/B^k(X)~$
denote the $k$-th reduced cohomology group of $X$ with real
coefficients. For each $k \geq -1$ endow $C^k(X)$ with the
standard inner product $(\phi,\psi)=\sum_{\sigma \in X(k)}
\phi(\sigma)\psi(\sigma)~~$ and the corresponding $L^2$ norm
$||\phi||=(\sum_{\sigma \in X(k)} \phi(\sigma)^2)^{1/2}$.
\\ Let $d_k^*:C^{k+1}(X) \rightarrow
C^k(X)$ denote the adjoint of $d_k$ with respect to these standard
inner products. The reduced $k$-Laplacian of $X$ is the mapping
$$\Delta_k=d_{k-1}d_{k-1}^*+d_k^*d_k : C^k(X) \rightarrow
C^k(X)~~.$$ Note that if $G$ denotes the $1$-skeleton of $X$ and
$J$ is the $V \times V$ all ones matrix, then the matrix $J+L_G$
represents $\Delta_0$ with respect to the standard basis. In
particular, the minimal eigenvalue of $\Delta_0$ equals
$\lambda_2(G)$.
\\
The space of harmonic $k$-cochains $\tilde{\cal H}^k(X) = \ker
\Delta_k$  consists of all $\phi \in C^k(X)$ such that both
$d_k\phi$ and $d_{k-1}^*\phi$ are zero. The simplicial version of
Hodge Theorem is the following well-known
\begin{proposition}
\label{hte} $\tilde{\cal H}^k(X) \cong \th^k(X)~$ for $k \geq 0$.
\end{proposition}
In particular, $\th^k(X)=0$ iff the minimal eigenvalue of
$\Delta_k$ is positive.

\section{Eigenvalues of Higher Laplacians}
\label{s:eandc}

Let $X=X(G)$ be the flag complex of a graph $G=(V,E)$ on $|V|=n$
vertices. For an $i$-simplex $\eta \in X$ let $\deg(\eta)$ denote
the number of $(i+1)$-simplices in $X$ which contain $\eta$. The
{\it link} of a simplex $\sigma \in X$ is the complex
$$\lk(\sigma)=\{\tau \in X ~:~\sigma \cup \tau \in X~,~ \sigma
\cap \tau = \emptyset~\}~.$$ For two ordered simplices $\sigma \in
X~,~\tau \in \lk(\sigma)$ let $\sigma\tau$ denote their ordered
union.
\begin{claim}
\label{dkpnorm} For $\phi \in C^k(X)$
$$
||d_k \phi||^2 =\sum_{\sigma \in X(k)} \deg(\sigma)
\phi(\sigma)^2 - 2 \sum_{\eta \in X(k-1)}\sum_{vw \in \lk(\eta)}
\phi(v \eta) \phi(w \eta)~~.
$$
\end{claim}
{\bf Proof:} Recall that for $\tau \in X(k+1)$ we denoted by
$\tau_i$ the ordered $k$-simplex obtained by removing the $i$-th
vertex of $\tau$. Thus $$||d_k \phi||^2=\sum_{\tau \in X(k+1)} d_k
\phi(\tau)^2=\sum_{\tau \in X(k+1)} \sum_{i=0}^{k+1} (-1)^i
\phi(\tau_i) \sum_{j=0}^{k+1} (-1)^j \phi(\tau_j)= $$ $$\sum_{\tau
\in X(k+1)} \sum_{i=0}^{k+1} \phi(\tau_i)^2 + \sum_{\tau\in
X(k+1)} \sum_{i \neq j} (-1)^{i+j} \phi(\tau_i) \phi(\tau_j)=$$ $$
\sum_{\sigma \in X(k)} \deg(\sigma) \phi(\sigma)^2 - 2 \sum_{\eta
\in X(k-1)}\sum_{vw \in \lk(\eta)} \phi(v \eta) \phi(w \eta)~~. $$
{\enp} For $\phi \in C^k(X)$ and a vertex $u \in V$ define $\phu
\in C^{k-1}(X)$ by $$\phu(\tau)= \left\{
\begin{array}{ll}
        \phi(u\tau)  & u \in \lk(\tau) \\
         0   & {\rm otherwise}
\end{array}
\right.~~ $$
\begin{claim}
\label{dphun}
For $\phi \in C^k(X)$ 
 $$\sum_{u \in V} ||d_{k-1}  \phu
||^2=$$
$$\sum_{\sigma \in X(k)} (\sum_{\tau \in \sigma(k-1)}
\deg(\tau)) \phi(\sigma)^2 - 2k\sum_{\tau \in X(k-1)} \sum_{vw \in
\lk(\tau)}\phi(v \tau)\phi(w \tau)~.
$$
\end{claim}
{\bf Proof:} Applying Claim \ref{dkpnorm} with $\phu \in
C^{k-1}(X)$ we obtain $$||d_{k-1}  \phu ||^2= \sum_{\tau \in
X(k-1)} \deg(\tau) \phu(\tau)^2-2 \sum_{\eta \in X(k-2)}\sum_{vw
\in \lk(\eta)} \phu(v \eta)\phu(w \eta)~.$$ Hence $$\sum_{u \in V}
||d_{k-1} \phu ||^2=$$ $$\sum_{u \in V} \sum_{\tau \in X(k-1)}
\deg(\tau) \phu(\tau)^2-2 \sum_{u \in V}\sum_{\eta \in
X(k-2)}\sum_{vw \in \lk(\eta)} \phu(v \eta)\phu(w \eta)=$$ $$
\sum_{\sigma \in X(k)} (\sum_{\tau \in \sigma(k-1)} \deg(\tau))
\phi(\sigma)^2 - 2\sum_{\eta \in X(k-2)}\sum_{vw \in \lk(\eta)}
\sum_{u \in \lk(v \eta) \cap \lk(w \eta)} \phi(vu \eta)\phi(wu
\eta)=$$ $$\sum_{\sigma \in X(k)} (\sum_{\tau \in \sigma(k-1)}
\deg(\tau)) \phi(\sigma)^2 - 2k \sum_{\tau \in X(k-1)}\sum_{vw \in
\lk(\tau)} \phi(v \tau)\phi(w \tau)~.$$ 
The last equality follows from
the fact that since $X$ is a flag complex, if $\eta \in
X(k-2)$~,~$vw \in \lk(\eta)$ and $u \in \lk(v \eta) \cap \lk(w
\eta)~$, then $vw \in \lk(u\eta)$.
 {\enp} Claims
\ref{dkpnorm} and \ref{dphun} imply $$
 k(||d_k \phi||^2 - \sum_{\sigma \in X(k)}
\deg(\sigma) \phi(\sigma)^2)=$$
\begin{equation}
\label{onee} \sum_{u \in V} ||d_{k-1}\phu||^2- \sum_{\sigma \in
X(k)}( \sum_{\tau \in \sigma(k-1)}\deg(\tau)) \phi(\sigma)^2~~.
\end{equation}
\begin{claim}
\label{phtn}
For $\phi \in C^k(X)$
\begin{equation}
\label{twoe}
 \sum_{u \in V} ||d_{k-2}^* \phu||^2= k ||\dkos
\phi||^2~~.
\end{equation}
\end{claim}
{\bf Proof:} For $\tau \in X(k-1)$ $$\dkos\phi(\tau)= \sum_{v \in
\lk(\tau)} \phi(v \tau)~.$$ Therefore $$ ||\dkos
\phi||^2=\sum_{\tau \in X(k-1)}\dkos \phi(\tau)^2=$$
\begin{equation}
\label{midt}
\sum_{\tau
\in X(k-1)}(\sum_{v \in \lk(\tau)} \phi(v \tau)) (\sum_{w \in
\lk(\tau)} \phi(w \tau))= \sum_{\tau \in X(k-1)}\sum_{(v,w) \in
\lk(\tau)^2} \phi(v \tau)\phi(w \tau)~.
\end{equation}
Substituting $\phu$ in (\ref{midt}) we obtain
 $$\sum_{u \in V} ||d_{k-2}^*\phu||^2=
\sum_{u \in V}\sum_{\eta \in X(k-2)}\sum_{(v,w) \in \lk(\eta)^2}
\phu(v \eta)\phu(w \eta)~=$$ $$ \sum_{\eta \in X(k-2)}\sum_{u \in
\lk(\eta)} \sum_{(v,w) \in \lk(u\eta)^2} \phi(v u\eta)\phi(w u
\eta)~= $$ $$ k\sum_{\tau \in X(k-1)}\sum_{(v,w) \in \lk(\tau)^2}
\phi(v \tau)\phi(w \tau)=k||\dkos \phi||^2~.$$ {\enp} Let $\phi
\in C^k(X)$. Summing (\ref{onee}) and (\ref{twoe}) we obtain the
following key identity: $$k(\Delta_k \phi,\phi)=$$
\begin{equation}
\label{thre}\sum_{u \in V}(\Delta_{k-1} \phu,\phu) - \sum_{\sigma
\in X(k)} (\sum_{\tau \in \sigma(k-1)} \deg(\tau)-k \deg(\sigma))
\phi(\sigma)^2~~.
\end{equation}
To estimate the righthand side of (\ref{thre}) we need the
following
\begin{claim}
\label{comb} For $\sigma \in X(k)$
\begin{equation}
\label{combi} \sum_{\tau \in \sigma(k-1)} \deg(\tau)-k
\deg(\sigma) \leq n~.
\end{equation}
\end{claim}
{\bf Proof:} Recall that $N(v)$ is the set of neighbors of $v$ in
$G$. Let $\sigma=[v_0,\ldots,v_k]$ then for any $I \subset
\{0,\ldots,k\}$ $$\deg([v_i:i \in I])=|\bigcap_{i \in I}
N(v_i)|~~.$$ Therefore
\begin{equation}
\label{clear} \sum_{\tau \in \sigma(k-1)} \deg(\tau)-k
\deg(\sigma)=\sum_{i=0}^k |\bigcap_{j \neq i} N(v_j)| -k
|\bigcap_{j=0}^k N(v_j)|~~.
\end{equation}
The Claim now follows since
each $v \in V$ is counted at most once on the righthand
side of (\ref{clear}). {\enp} {\bf Proof of
Theorem \ref{eigenv}:} Let $0 \neq \phi \in C^k(X)$ be an
eigenvector of $\Delta_k$ with eigenvalue $\mu_k(G)~$. By double
counting
\begin{equation}
\label{nor} \sum_{u \in V} ||\phu||^2=(k+1)||\phi||^2~~.
\end{equation}
Combining (\ref{thre}),(\ref{combi}) and (\ref{nor}) we obtain $$
k \mu_k(G) ||\phi||^2= k(\Delta_k \phi,\phi) \geq \sum_{u \in
V}(\Delta_{k-1} \phu,\phu) - n \sum_{\sigma \in X(k)}
\phi(\sigma)^2~~  \geq $$ $$\mu_{k-1}(G) \sum_{u \in V} ||\phu||^2
- n ||\phi||^2 = ((k+1)\mu_{k-1}(G)-n) ||\phi||^2~~. $$ {\enp} \
\\
\\ {\bf Proof of Theorem \ref{ei}:} Inequality (\ref{maineq})
implies by induction on $k$ that $\mu_k(G) \geq (k+1) \mu_0(G)
-kn$. Therefore, if $\mu_0(G)=\lambda_2(G)>\frac{kn}{k+1}$ then
$\mu_k(G)>0$ and $\th^{k}(X(G),\Rea)=0$ follows from the
simplicial Hodge Theorem. {\enp}

\section{Vector Domination and Homology}
\label{s:etag} Let $G=(V,E)$ be a graph with $|V|=n$. We first
reformulate Theorem \ref{ei} in terms of the independence complex
$\cig$.
\begin{theorem}
\label{indeta} $\eta(\cig) \geq \frac{n}{\lambda_n(G)}.$
\end{theorem}
{\bf Proof:} Let $\ell=\lceil\frac{n}{\lambda_n(G)}\rceil$. Since
$\lambda_n(G)=n-\lambda_2(\comg)$ it follows that
$\lambda_2(\comg)>\frac{\ell-2}{\ell-1}n$. Therefore by Theorem
\ref{ei}, $\th^i(\cig)=\th^i(X(\comg))=0$ for $i \leq \ell-2$.
Hence $\eta(\cig) \geq \ell$. {\enp} The proof of Theorem
\ref{etagegamma} depends on Theorem \ref{indeta} and the following
\begin{claim}
\label{lambdan} Let $P$ be a vector representation of $G=(V,E)$.
Then
$$
\lambda_n(G) \leq \max_{u \in V}~ P(u) \cdot
\sum_{v \in V} P(v)~~.
$$
\end{claim}
{\bf Proof:} Let $x=(x(v):v \in V)$ be a vector in $\Rea^V$. Then
$$x^T L_G x =\sum_{uv \in E} (x(u)-x(v))^2 \leq $$ $$\frac{1}{2}
\sum_{(u,v) \in V \times V} (x(u)-x(v))^2 P(u) \cdot P(v) =$$ $$
\sum_{u \in V}x(u)^2P(u) \cdot \sum_{v \in V} P(v)  ~-~ ||
\sum_{v\in V} x(v)P(v)||^2 \leq $$ $$ ||x||^2~ \max_{u \in V}~
P(u) \cdot \sum_{v \in V} P(v)~.$$ The Claim follows since
$\lambda_n(G)= \max~\{\frac{x^TL_Gx}{||x||^2}~:~0 \neq x \in
\Rea^V\}.$ {\enp} Let $\Int_+$ denote the positive integers and
let $\Rat_+$ denote the positive rationals. For a vector
$\va=(a(v):v \in V) \in \Int_+^V$ let $\gva$ denote the graph
obtained by replacing each $v \in V$ by an independent set of size
$a(v)$. Formally $V(\gva)=\{(v,i):v \in V~,~1 \leq i \leq a(v)\}$
and $\{(u,i),(v,j)\} \in E(\gva)$ if $\{u,v\} \in E$. The projection
$(v,i) \rightarrow v$ induces a homotopy equivalence between
$\ci(\gva)$ and $\ci(G)$. In particular $
\eta(\ci(\gva))=\eta(\ci(G))$.
\
\\
\\ {\bf Proof of Theorem \ref{etagegamma}:} Let $P$ be a
representation of $G$. By linear programming duality $$ |P|=\min\{
\valpha \cdot \vone~:~\valpha \geq 0~,~\valpha PP^T \ge
\vone~\}=$$ $$\max\{ \valpha \cdot \vone~:~\valpha \geq 0~,~\valpha
PP^T \leq \vone~\}=$$
$$
{\rm Sup}\{\valpha \cdot
\vone~:~\valpha \in \Rat_+^V~~,~~\valpha PP^T \leq \vone~\}~.
$$ Let $\valpha \in \Rat_+^V$ such that $\valpha PP^T \leq \vone$.
Write $\valpha= \frac{\va}{k}$ where $k \in \Int_+$ and
$\va=(a(v):v \in V) \in \Int_+^V$. Let $N=|V(\gva)|=\sum_{u \in
V}a(u)$. Consider the representation $Q$ of $\gva$ given by
$Q((u,i))=P(u)$ for $(u,i) \in V(\gva)$. By Claim \ref{lambdan}
$$\lambda_N(\gva) \leq \max_{(u,i) \in V(\gva)}~ Q((u,i)) \cdot
\sum_{(v,j) \in V(\gva)} Q((v,j))= $$ $$ \max_{u \in V}~P(u) \cdot
\sum_{v \in V} a(v)P(v) \leq k~~.$$ Hence by Theorem \ref{indeta}
$$\valpha \cdot \vone= \frac{1}{k}\sum_{v \in V} a(v) =
\frac{N}{k} \leq$$ $$ \frac{N}{\lambda_N(\gva)} \leq
\eta(\ci(\gva))=\eta(\ci(G))~~.$$ {\enp} \noindent {\bf Remarks:}
\vspace{5pt}
\\ 1. Let $C_n$ denote the $n$-cycle on the vertex
set $V=\{0,\ldots, n-1\}$.  For $n=3k$ define a representation $P$
of $C_{3k}$ by $$ P(\ell)= \left\{
\begin{array}{ll}
        e_{2j} & \ell=3j \\
        e_{2j}+e_{2j+1} & \ell=3j+1 \\
        e_{2j+1}+e_{2j+2} & \ell=3j+2
\end{array}
\right. $$ where $e_0,\ldots,e_{2k-1}$ are orthogonal unit vectors
 and the indices are cyclic modulo $2k$. Let $\alpha
\in \Rea^V$ be given by $~\alpha(\ell)=1~$ if $3$ divides $\ell$
and zero otherwise. Since $\alpha PP^T = \vone$, it follows by
linear programming duality that $\Gamma(C_{3k}) \geq \alpha \cdot
\vone=k$. On the other hand (see Claim 3.3 in \cite{M03})
$\eta(I(C_n))=\lfloor\frac{n+1}{3}\rfloor$. Therefore
$\eta(\ci(C_{3k}))=\Gamma(C_{3k})=k$. For $n=3k+1$ it can
similarly be shown that $\eta(\ci(C_{3k+1}))=\Gamma(C_{3k+1})=k$.
The case $n=3k-1$ is more involved and we only have the bounds
$k-\frac{1}{2} \leq \Gamma(C_{3k-1}) \leq \eta(\ci(C_{3k-1}))=k$.
Note that for cycles the bound $\eta(\ci(G)) \geq \gamma^*_s(G)$
is weaker since $\gamma^*_s(C_n)=\frac{n}{4}$.  \vspace{1pt} \\ 2.
It can be shown that for any graph $\Gamma(G) \geq {\rm
Sup}\{\gamma^*_s(\gva)~:~\va \in \Int_+^V\}. $ We do not know of
examples with strict inequality.

\section{A Hall Type Theorem for Fractional Width}
\label{s:hall}

Let $Z$ be a simplicial complex on the vertex set $W$ and let
$\bigcup_{i=1}^m W_i$ be a partition of $W$. A simplex $\tau \in
Z$ is {\it colorful} if $|\tau \cap W_i|=1$ for all $1 \leq i \leq
m$. For $W' \subset W$ let $Z[W']$ denote the induced subcomplex
on $W'$. The following Hall's type sufficient condition for the
existence of colorful simplices appears in \cite{AH00} and in
\cite{M01}.
\begin{proposition}
\label{hom} If for all $\emptyset \neq I \subset [m]$
$$\eta(Z[\bigcup_{i \in I} W_i]) \geq |I|$$ then $Z$ contains a
colorful simplex.
\end{proposition}
Let $G$ be a graph on the vertex set $W$ with a partition
$W=\bigcup_{i=1}^m W_i$. A set $S \subset W$ is {\it colorful} if
$S \cap W_i \neq \emptyset$ for all $1 \leq i \leq m$. The induced
subgraph on $W' \subset W $ is denoted by $G[W']$.  Combining
Theorem \ref{etagegamma} and Proposition \ref{hom} we obtain the
following
\begin{theorem}
\label{gammahall} If  $\Gamma(G[\bigcup_{i\in I} W_i]) > |I|-1$
for all $\emptyset \neq I \subset [m]$ then $G$ contains a
colorful independent set.
\end{theorem}
Let $\cf \subset 2^V$ be a hypergraph, possibly with multiple
edges. The {\it line graph} $G_{\cf}=(W,E)$ associated with $\cf$
has vertex set $W=\cf$ and edge set $E$ consisting of all
$\{F,F'\} \subset \cf$ such that $F \cap F' \neq \emptyset$. A
matching in $\cf$ corresponds to an independent set in $G_{\cf}$.
For each $F \in \cf$ let $P(F) \in \Rea^V$ denote the incidence
vector of $F$. $P$ is clearly a vector representation of $G_{\cf}$
and satisfies $|P|=w^*(\cf)$. Thus $\Gamma(G_{\cf}) \geq
w^*(\cf)$. \ \\
\\ {\bf Proof of
Theorem \ref{wstar}:} Let $\cf$ denote the disjoint union of the
$\cf_i$'s, and consider the graph $G_{\cf}=(W,E)$ with the partition
$W=\cup_{i=1}^m W_i$ where $W_i=\cf_i$. Then for any $\emptyset
\neq I \subset [m]$ $$\Gamma(G_{\cf}[\cup_{i \in I}
W_i])=\Gamma(G_{\cup_{i \in I} \cf_i}) \geq $$ $$w^*(\cup_{i \in
I} \cf_i) > |I|-1~~.$$ Theorem \ref{gammahall} implies that
$G_{\cf}$ contains a colorful independent set, hence $\cfi$
contains an SDR. {\enp}

\end{document}